\documentclass{article}
\usepackage[centertags]{amsmath}
\usepackage{amsfonts}
\usepackage{amssymb}
\usepackage{amsthm}
\usepackage{newlfont}

\font\mostafa=cmr12
\font\ramadan=cmr10
\title{\bf On the Approximate Solution of Integral Equations with Logarithmic Kernels Using the Third Kind of Chebyshev Polynomials}
\author{M. R. A. Sakran\footnote{Corresponding author. E-mail:
$mrs01@fayoum.edu.eg$}~\\Department of
Mathematics, Faculty of Science, Fayoum University, Fayoum,
Egypt.}
\date{}
\begin{document}
\maketitle
\begin{flushleft}
{\bf Abstract}
\end{flushleft}

{\ramadan \ \ \ \ \  An expansion procedure using third kind Chebyshev polynomials as base functions is suggested for solving second type Volterra integral equations with logarithmic kernels. The algorithm's convergence is studied and some illustrative examples are presented to show the method's efficiency and reliability,  comparisons with other methods in the literature are made.\\


{\bf Mathematics Subject Classification (MSC)}: 45D05, 45E10, 65R20, 33C47, 65L20, 65L70.
\section{Introduction}
\mostafa \ \ \ \ \ \ \ \ \ Numerous fields of numerical analysis frequently employ the Chebyshev polynomials. In [9,10,11,17], a thorough explanation of the uses and properties of Chebyshev polynomials is given.

Let us consider the following second kind Volterra integral equation (VIE) with logarithmic kernel:
\begin{equation}
y(x)=f(x)+\int_{0}^{x}\ln(x-s)k(x,s,y(s))ds,
\end{equation}
We assume that (1) has a unique solution $y{\in}C[0,X]$ that can be determined, see [4,5,6]. These equations are usually difficult to solve analytically, thus numerical techniques must be used to arrive at an accurate approximation.

Equation (1) has been solved using a variety of numerical techniques; for example, see [1,2,4,7,13,14,15,16,18,19,20,22,23].

This article aims to create a computational approach that represents the unknown function by a finite expansion of Chebyshev polynomials of the third kind to solve equation (1).

In the next section, the proposed strategy is derived. The third section discusses the suggested method's convergence analysis. The main theoretical results obtained in the third section are demonstrated in the fourth section by some numerical tests. We end the paper with some closing remarks.
\section{Derivation of the algorithm}

\ \  Let $\mathbb{R}=(-\infty,\infty), \ \mathbb{N}=\{1,2,3,...\}$ \ and $\mathbb{N}_{0}=\mathbb{N}\bigcup\{0\}$. We establish a mesh, ${P_{N}}=\{0=x_{0}<x_{1}<...<x_{N}=X\}$, for $N\in\mathbb{N}$ with a step length $h_{i}=x_{i+1}-x_{i}, i=0(1)N-1(i=0,1,2,...,N-1)$, for the interval $J=[0,X]$. The points $x_{ij}=x_{i}+\beta_{j}h_{i}$ split each subinterval $[x_{i},x_{i+1}]$ where\\
\begin{equation}
    \beta_{j}=\frac{1}{2}(1-\cos(\frac{j\pi}{\nu})), \ \ j=0(1)\nu.
\end{equation}
Assume the following approximation for $x \in [x_{i},x_{i+1}]$
\begin{equation}
  y(x)\approx u_{i+1}(x)=\sum_{r=0}^{ \ \nu }a_r^{(i+1)}V_{r}[\frac{2}{h_{i}}(x-x_{i})-1],
\end{equation}
where [21]
\begin{equation}
a_r^{(i+1)}=\frac{1}{\nu}\sum_{k=0}^{ \ \nu \ \ _{''}}(1-\cos(\frac{k\pi}{\nu}))V_{r}[-\cos(\frac{k\pi}{\nu})]u_{i+1,k}, \ \ r=0(1)\nu.
\end{equation}
Where $V_{r}[x]$ is the third kind Chebyshev polynomial of degree $r$ [9,12,17], a summation sign with double primes indicates a sum with first and last terms halved, and $u_{i+1,k}= u_{i+1}(x_{ik})$.

Integrating (3), multiplied by $\ln(x_{ij}-x)$, from $x_{i}$ to $x_{ij}$ relative to $x$, and then setting $x=x_{i}+\frac{1}{2}h_{i}(1-\cos(\frac{j\pi}{\nu})-t)$ results in
\begin{equation}
\int_{x_{i}}^{x_{ij}}\ln(x_{ij}-x)u_{i+1}(x)dx=\frac{h_{i}}{2}
\sum_{r=0}^{ \ \nu }(-1)^{r}a_r^{(i+1)}I_{r,j}(\cos(\frac{j\pi}{\nu})),
\end{equation}
where
\begin{equation}
I_{r,j}(\eta)=\int_{0}^{1-\eta}\ln(\frac{h_{i}}{2}t) W_r[t+\eta]dt.
\end{equation}
Here $W_{r}[x]$ is the r-th Chebyshev polynomial of the fourth kind [9,12,17].\\
To determine $I_{r,j}(\eta)$: Let,
\begin{eqnarray}
  G &=& \int_{0}^{1-\eta}t(t+\eta-1)\ln(\frac{h_{i}}{2}t) W_r[t+\eta]dt \\
 \nonumber  &=& \int_{0}^{1-\eta}[(t+\eta)-\eta][(t+\eta)-1]\ln(\frac{h_{i}}{2}t)W_r[t+\eta]dt \\
 \nonumber   &=& \int_{0}^{1-\eta}\ln(\frac{h_{i}}{2}t)(t+\eta)^{2}W_r[t+\eta]dt-(1+\eta)\int_{0}^{1-\eta}\ln(\frac{h_{i}}{2}t)
(t+\eta)W_r[t+\eta]dt \\
 \nonumber  &+& \eta\int_{0}^{1-\eta}\ln(\frac{h_{i}}{2}t)W_r[t+\eta]dt,
\end{eqnarray}
then,
\begin{equation}
G=\frac{1}{4}(I_{r+2,j}(\eta)+2I_{r,j}(\eta)+I_{r-2,j}(\eta))-\frac{1+\eta}{2}(I_{r+1,j}(\eta)+I_{r-1,j}(\eta))+\eta
I_{r,j}(\eta).
\end{equation}
Integration by parts reduces (7) as
\begin{eqnarray*}
  G &=& \frac{-1}{2(r+1)}\int_{0}^{1-\eta}\ln(\frac{h_{i}}{2}t)(t+\eta-1)W_{r+1}[t+\eta]dt \\
    &-& \frac{1}{2r(r+1)}\int_{0}^{1-\eta}\ln(\frac{h_{i}}{2}t)(t+\eta-1)W_{r}[t+\eta]dt \\
    &+& \frac{1}{2r}\int_{0}^{1-\eta}\ln(\frac{h_{i}}{2}t)(t+\eta-1)W_{r-1}[t+\eta]dt \\
    &-& \frac{1}{2(r+1)}\int_{0}^{1-\eta}\ln(\frac{h_{i}}{2}t)tW_{r+1}[t+\eta]dt \\
    &-& \frac{1}{2r(r+1)}\int_{0}^{1-\eta}\ln(\frac{h_{i}}{2}t)tW_{r}[t+\eta]dt \\
    &+& \frac{1}{2r}\int_{0}^{1-\eta}\ln(\frac{h_{i}}{2}t)tW_{r-1}[t+\eta]dt \\
    &-& \frac{1}{2(r+1)}\int_{0}^{1-\eta}(t+\eta-1)W_{r+1}[t+\eta]dt \\
    &-& \frac{1}{2r(r+1)}\int_{0}^{1-\eta}(t+\eta-1)W_{r}[t+\eta]dt \\
    &+& \frac{1}{2r}\int_{0}^{1-\eta}(t+\eta-1)W_{r-1}[t+\eta]dt,
\end{eqnarray*}
i.e.,
\begin{eqnarray*}
G &=& \frac{-1}{4(r+1)}(I_{r+2,j}(\eta)+I_{r,j}(\eta))+\frac{1}{2(r+1)}I_{r+1,j}(\eta) \\
&-& \frac{1}{4r(r+1)}(I_{r+1,j}(\eta)+I_{r-1,j}(\eta))+\frac{1}{2r(r+1)}I_{r,j}(\eta) \\
&+& \frac{1}{4r}(I_{r,j}(\eta)+I_{r-2,j}(\eta))-\frac{1}{2r}I_{r-1,j}(\eta) \\
&-& \frac{1}{4(r+1)}(I_{r+2,j}(\eta)+I_{r,j}(\eta))+\frac{\eta}{2(r+1)}I_{r+1,j}(\eta) \\
&-& \frac{1}{4r(r+1)}(I_{r+1,j}(\eta)+I_{r-1,j}(\eta))+\frac{\eta}{2r(r+1)}I_{r,j}(\eta) \\
\end{eqnarray*}
\begin{eqnarray}
\nonumber    &+& \frac{1}{4r}(I_{r,j}(\eta)+I_{r-2,j}(\eta))-\frac{\eta}{2r}I_{r-1,j}(\eta) \\
\nonumber    &-& \frac{1}{4(r+1)}(Q_{r+2,j}(\eta)+Q_{r,j}(\eta))+\frac{1}{2(r+1)}Q_{r+1,j}(\eta) \\
\nonumber    &-& \frac{1}{4r(r+1)}(Q_{r+1,j}(\eta)+Q_{r-1,j}(\eta))+\frac{1}{2r(r+1)}Q_{r,j}(\eta) \\
    &+& \frac{1}{4r}(Q_{r,j}(\eta)+Q_{r-2,j}(\eta))-\frac{1}{2r}Q_{r-1,j}(\eta)
\end{eqnarray}
The following recurrence relation, which can be constructed by equating (8) and (9) and then multiplying by $4r(r+1)$, may be used to get the value of the definite integral (6).
\begin{eqnarray}
\nonumber  && r(r+3)I_{r+2,j}(\eta)+2[1-(1+ \eta)r(r+2)]I_{r+1,j}(\eta)+2[(1+2\eta)r(r+1)-\eta-2]I_{r,j}(\eta) \\
\nonumber     &+& 2[1-(1+\eta)(r^{2}-1)]I_{r-1,j}(\eta)+(r+1)(r-2)I_{r-2,j}(\eta)=-rQ_{r+2,j}(\eta)+(2r-1)Q_{r+1,j}(\eta) \\
    &+& 3Q_{r,j}(\eta)-(2r+3)Q_{r-1,j}(\eta)+(r+1)Q_{r-2,j}(\eta), \ \ \ \ \ r=2(1)\nu-2, \ j=0(1)\nu,
\end{eqnarray}
where $Q_{r,j}(\eta)=\int_{0}^{1-\eta}W_r[t+\eta]dt$\\
i.e.
\begin{equation}
2r(r+1)Q_{r,j}(\eta)=2(2r+1)-rW_{r+1}[\eta]-W_{r}[\eta]+(r+1)W_{r-1}[\eta], \ \ \
\  r=2(1)\nu, \ j=0(1)\nu,
\end{equation}
\begin{equation}
Q_{0,j}(\eta)=1-\eta, \ \ \ Q_{1,j}(\eta)=(1-\eta)(2+\eta).
\end{equation}
Initial values are as follows:
\begin{eqnarray}
  I_{0,j}(\eta) &=& 2\beta_{j}[\ln(h_{i}\beta_{j})-1] \\
  I_{1,j}(\eta) &=& \beta_{j}[2(2+\eta)\ln(h_{i}\beta_{j})-3(1+\eta)] \\
  I_{2,j}(\eta) &=& \frac{1}{9}\beta_{j}[6(4\eta^{2}+7\eta+4)\ln(h_{i}\beta_{j})-(44\eta^{2}+47\eta-1)] \\
  I_{3,j}(\eta) &=& \frac{1}{9}\beta_{j}[6(6\eta^{3}+10\eta^{2}+4\eta+1)\ln(h_{i}\beta_{j})-(75\eta^{3}+83\eta^{2}-13\eta-19)].
\end{eqnarray}
By inserting the formulae of $I_{r,j}(\cos\frac{j\pi}{\nu})$ in (5) and using (4), yields after some calculations
\begin{equation}
\bigl[\int_{x_{i}}^{x_{ij}}\ln(x_{ij}-x)u_{i+1}(x)dx\bigr]=h_{i}C[u_{i+1,j}],
\end{equation}
where the matrix $C$'s elements are described by
\begin{equation}
c_{jk}=\frac{1-\cos(\frac{k\pi}{\nu})}{2\nu(1+\delta_{k0}+\delta_{k\nu})}\sum_{r=0}^{ \
\nu}(-1)^{r}I_{r}(\cos(\frac{j\pi}{\nu}))V_{r}[-\cos(\frac{k\pi}{\nu})],
\ \ \  j=0(1)\nu, \ k=0(1)\nu,
\end{equation}
where $\delta_{ij}$ is the Kronecker delta.

Using (17), the numerical solution to (1) may be determined as follows:\\
Put $x=x_{ij}\in[{x_{i}},{x_{i+1}}]$ in (1) to obtain
\begin{eqnarray}
\nonumber  y(x_{ij}) &=& f(x_{ij})+\int_{0}^{x_{ij}}\ln(x_{ij}-s)k(x_{ij},s,y(s))ds \\
\nonumber    &=& f(x_{ij})+\sum_{l=0}^{i-1}\int_{x_{l}}^{x_{l+1}}\ln(x_{ij}-s)k(x_{ij},s,y(s))ds+\int_{x_{i}}^{x_{ij}}\ln(x_{ij}-s)k(x_{ij},s,y(s))ds,\\
&& \ \ \ \ \ \ \ \ \ \ \ \ \ \ \ \ \ \ \ \ \ \ \ \ \ \ \ \ \ \ \ \ \ \ \ \ \ \ \ \ \ \ \ \ \ \ \ \ \ \ i=0(1)N-1, \  j=1(1)\nu,
\end{eqnarray}
Using (17) in (19), the result is
\begin{equation}
u_{i+1,j}=f(x_{ij})+\sum_{l=0}^{i-1}h_{l}\sum_{k=0}^{\nu}b_{\nu
k}\ln(x_{ij}-x_{lk})k(x_{ij},x_{lk},u_{l+1,k})+h_{i}\sum_{k=0}^{\nu}c_{jk}k(x_{ij},x_{ik},u_{i+1,k}),
\end{equation}
where [21]
\begin{equation}
    b_{{\nu}k}=\frac{1-\cos(\frac{k\pi}{\nu})}{2\nu(1+\delta_{k0}+\delta_{k{\nu}})}\{2+\sum_{r=1}^{ \ \nu }\frac{(-1)^{r}(2r+1)-1}{r(r+1)}V_{r}[-\cos(\frac{k\pi}{\nu})]\},
\end{equation}
\begin{equation*}
   u_{1,0}=f(0), \ u_{i+1,0}=u_{i,\nu}, \ i=0(1)N-1,  \ j=1(1)\nu.
\end{equation*}

The Gaussian elimination approach for a linear system and the Newton iterative method for a nonlinear system of equations may both be used to find the $\nu$-unknowns $u_{i+1,j}$, $j=1(1)\nu$ of (20) in each subinterval $[x_{i},x_{i+1}]$, $i=0(1)N-1$.\\
\underline
{\bf Remark: The method's cost}

As in [21], we can conclude that the proposed method costs $O((\nu+1)^{3})$ operations for each subinterval, i.e., $N \times O((\nu+1)^{3})$ operations on $J$. So, the cost is effectively controlled since $N$ and $\nu$ are small.


\section{Convergence Analysis}
\mostafa
 \ \ \ Results for $T< e$ and the linear version of (1) are derived here, namely

\begin{equation}
y(x)=f(x)+\int_{0}^{x}\ln(x-s)K(x,s)y(s)ds.
\end{equation}
By rewriting (22) as the following for $x\in[x_{i},x_{i+1}]$
\begin{equation}
y(x)=f(x)+\sum_{l=0}^{i-1}\int_{x_{l}}^{x_{l+1}}\ln(x-s)K(x,s)y(s)ds+\int_{x_{i}}^{x}\ln(x-s)K(x,s)y(s)ds.
\end{equation}
Put $x=x_{ij}\in[{x_{i}},{x_{i+1}}]$, we get
\begin{equation}
y(x_{ij})=f(x_{ij})+\sum_{l=0}^{i-1}\int_{x_{l}}^{x_{l+1}}\ln(x_{ij}-s)K(x_{ij},s)y(s)ds+\int_{x_{i}}^{x_{ij}}\ln(x_{ij}-s)K(x_{ij},s)y(s)ds.
\end{equation}
Equation (24) is satisfied approximately by a numerical approximation $u_{i+1}$ to the exact solution $y$ of
(22), so
\begin{equation*}
e(x_{ij})=\sum_{l=0}^{i-1}\int_{x_{l}}^{x_{l+1}}\ln(x_{ij}-s)K(x_{ij},s)e_{l}(s)ds+\int_{x_{i}}^{x_{ij}}\ln(x_{ij}-s)K(x_{ij},s)y(s)ds-\int_{x_{i}}^{x_{ij}}\ln(x_{ij}-s)K(x_{ij},s)u_{i+1}(s)ds,
\end{equation*}
since $X\leq e$, then ${\ln}(x_{ij}-s)\leq 1 \ \ \ \
{\forall}s\in[{x_{l}},{x_{l+1}}], l=0(1)i-1$, then
\[e(x_{ij})\leq \sum_{l=0}^{i-1}h_{l}\int_{0}^{1}K(x_{ij},x_{l}+vh_{l})e_{l}(x_{l}+vh_{l})dv+\int_{x_{i}}^{x_{ij}}\ln(x_{ij}-s)K(x_{ij},s)y(s)ds-\int_{x_{i}}^{x_{ij}}\ln(x_{ij}-s)K(x_{ij},s)u_{i+1}(s)ds,\]
where
\begin{eqnarray*}
  e(x_{ij}) &=& y(x_{ij})-u_{i+1}(x_{ij}), i=0(1)N-1, j=0(1)\nu,\\
  e_{l}(x_{l}+vh_{l}) &=& y(x_{l}+vh_{l})-u_{i+1}(x_{l}+vh_{l}), l=0(1)i-1.
\end{eqnarray*}
If \ \ \ $k_{0}:={\max}\{{\mid}K(x,s)\mid:(x,s){\in}\mathfrak{F} ( where \ \mathfrak{F}: \{(x,s): 0{\leq}s{\leq}x{\leq}X\})\}$,\\
then
\begin{equation}
{\mid}e(x_{ij}){\mid}\leq
k_{0}\sum_{l=0}^{i-1}h_{l}\int_{0}^{1}{\mid}e_{l}(x_{l}+vh_{l}){\mid}dv+k_{0}{\mid}\rho_{ij}{\mid},
\end{equation}
where
\begin{equation}
\rho_{ij}=\int_{x_{i}}^{x_{ij}}\ln(x_{ij}-s)y(s)ds-\int_{x_{i}}^{x_{ij}}\ln(x_{ij}-s)u_{i+1}(s)ds.
\end{equation}

In the following, we will deal with three types of mesh sequences, which are uniform, quasi-uniform, and graded meshes. For any one of them, we find that $h=O(N^{-1})$ for each compact interval $J$, where $h=\max_{i}h_{i}, \  0{\leq}i{\leq}N-1$, is the mesh diameter for a partition ${P_{N}}$ of $J$ see [3].
\newtheorem{lemma}{Lemma}[section]
\begin{lemma}
\mostafa If $x_{i},i=0(1)N$, is a sequence of real numbers with
\begin{equation}
\mid x_{i}\mid\leq hM \sum_{j=0}^{i-1}\mid x_{j}\mid+\delta, \ \
i=1(1)N,
\end{equation}
where $M>0$ and usually independent of $h$ and $\delta>0$ then
\begin{equation}
\mid x_{i}\mid\leq (hM \mid x_{0}\mid+\delta)e^{Mih}, \ \ i=1(1)N.
\end{equation}
\end{lemma}
\mostafa Proof: see [8].
\begin{lemma}
For $j\in\mathbb{N}$, we have
\begin{equation}
\sum_{k=0}^{ \ \nu \ \
_{''}}\beta_{k}^{j}V_{r}[-cos(\frac{k\pi}{\nu})]=\left\{\begin{array}{ll}
                \frac{1}{2}&\mbox{if $\nu$=1$ \ \ \forall$r,j},\\
                0&\mbox{if r$\geq$ j,r=j$\not=$$\nu$},\\
                 \frac{\nu}{2^{2j-1}}[(_{j-r-1}^{2j-1})+(_{j-r-2}^{2j-1})\delta_{r,\nu-1}+(_{j-r}^{2j-1})\delta_{r\nu}]&\mbox{if r${<}$ j,r=j=$\nu$},
                 \end{array}
        \right.
\end{equation}
where $(_{r}^{n})=\frac{n!}{r!(n-r)!}$.
\end{lemma}
Proof: see [21].
\begin{lemma}
\mostafa If the functions $f$ and $K$ in (22) are such that
$f{\in}C^{\nu-1}(J)$ and $K{\in}C^{\nu-1}(\mathfrak{F})$,and
$y^{(\aleph)}(x_{i})=u_{i+1}^{(\aleph)}(x_{i}) \ \forall \aleph, \aleph \in\mathbb{N}_{0}$, then
\begin{equation}
{\mid}\rho_{ij}\mid \leq \mathcal{P} h^{\nu},
\end{equation}
for any mesh sequence, where $\mathcal{P}$ is a generic constant depends on
the $\{\beta_{j}\}$.
\end{lemma}
Proof:\\
By (26)
\begin{eqnarray*}
\rho_{ij}&=&\int_{x_{i}}^{x_{ij}}\ln(x_{ij}-s)y(s)ds-\int_{x_{i}}^{x_{ij}}\ln(x_{ij}-s)u_{i+1}(s)ds\\
&=&\int_{x_{i}}^{x_{ij}}\ln(x_{ij}-s)y(s)ds-h_{i}\sum_{k=0}^{\nu}c_{jk}u_{i+1,k}\\
&=&h_{i}\{\int_{0}^{\beta_{j}}\ln(h_{i}\Psi)y(x_{i}+(\beta_{j}-\Psi)h_{i})d{\Psi}-\sum_{k=0}^{\nu}c_{jk}u_{i+1}(x_{i}+\beta_{k}h_{i})\}.
\end{eqnarray*}
Using the Taylor expansion results in
\begin{equation*}
\rho_{ij}=h_{i}\sum_{n=0}^{\nu-1}\frac{h_{i}^{n}}{n!}y^{(n)}(x_{i})[\int_{0}^{\beta_{j}}\ln(h_{i}\Psi)(\beta_{j}-\Psi)^{n}d{\Psi}-\sum_{k=0}^{\nu}\beta_{k}^{n}c_{jk}],
\end{equation*}
then
\begin{equation}
\rho_{ij}=h_{i}\sum_{n=0}^{\nu-1}\frac{h_{i}^{n}}{n!}y^{(n)}(x_{i})[\Omega(n)-\sum_{k=0}^{\nu}\beta_{k}^{n}c_{jk}],
\end{equation}
where
\begin{equation}
\Omega(n)=\beta_{j}^{n+1}\sum_{\ell=0}^{n}\frac{(-1)^{\ell}}{(\ell+1)^{2}}[(\ell+1)\ln(h_{i}\beta_{j})-1].
\end{equation}
Now, we look for the first value of $n$ that outcomes in
\begin{equation}
\sum_{k=0}^{\nu}\beta_{k}^{n}c_{jk}\not=\Omega(n).
\end{equation}
If this value is $\mu$, then ${\mid}\rho_{ij}\mid \leq
\mathcal{P}h^{\mu+1}$,
where $\mathcal{P}$ is a generic constant depending on the $\{\beta_{j}\}$.\\
Now, we have
\begin{equation*}
\sum_{k=0}^{\nu}\beta_{k}^{n}c_{jk}=\frac{1}{\nu}\sum_{r=0}^{\nu}(-1)^{r}I_{r}(cos(\frac{j\pi}{\nu}))\sum_{k=0}^{ \ \nu \ \
_{''}}\beta_{k}^{n+1}V_{r}[-cos(\frac{k\pi}{\nu})].
\end{equation*}
Lemma 3.2 gives:\\
$\ast$ \ If $\nu=1$ results in
\begin{equation}
\sum_{k=0}^{\nu}\beta_{k}^{n}c_{jk}\not=\Omega(n), \ \forall n.
\end{equation}
$\ast$ \ If $r {\geq} n+1, r=n+1\not=\nu$, We get
\begin{equation}
\sum_{k=0}^{\nu}\beta_{k}^{n}c_{jk}=0\not=\Omega(n), \ \forall n.
\end{equation}
$\ast$ \ If $r < n+1$ results in\\
$\bullet$ \ For $n\leq\nu-2$,  yields
\begin{eqnarray*}
\sum_{k=0}^{\nu}\beta_{k}^{n}c_{jk}&=&\frac{1}{2^{2n+1}}\int_{0}^{1-\eta}ln(\frac{h_{i}}{2}x)
\sum_{r=0}^{n}(^{2n+1}_{n-r})V_{r}[-(x+\eta)]dx\\
&=&\frac{1}{2^{n+1}}\sum_{\jmath=0}^{n}(^{n}_{\jmath})(-1)^{\jmath}(1-\eta)^{n-\jmath}\int_{0}^{1-\eta}x^{\jmath}ln(\frac{h_{i}}{2}x)dx,
\end{eqnarray*}
then
\begin{equation}
\sum_{k=0}^{\nu}\beta_{k}^{n}c_{jk}=\Omega(n), \ \ n\leq\nu-2.
\end{equation}
$\bullet$ For the other cases: \ If $n\geq\nu-1$ or $\ r=n+1=\nu$, we have\\
\begin{eqnarray*}
\sum_{k=0}^{\nu}\beta_{k}^{n}c_{jk}&=&\frac{1}{2^{2n+1}}\int_{0}^{1-\eta}ln(\frac{h_{i}}{2}x)
\sum_{r=0}^{\nu}[(^{2n+1}_{n-r})+(^{2n+1}_{n-r+1})\delta_{r\nu}+(^{2n+1}_{n-r-1})\delta_{r,\nu-1}]V_{r}[-(x+\eta)]dx\\
&=&\Omega(n)-\frac{1}{2^{2n+1}}\int_{0}^{1-\eta}ln(\frac{h_{i}}{2}x)\{\sum_{l=\nu+1}^{n}(^{2n+1}_{n-l})V_{l}[-(x+\eta)]\\
&-&(^{2n+1}_{n-\nu})V_{\nu-1}[(-(x+\eta)]-(^{2n+1}_{n-\nu+1})V_{\nu}[(-(x+\eta)]\}dx,
\end{eqnarray*}
then
\begin{equation}
\sum_{k=0}^{\nu}\beta_{k}^{n}c_{jk}\not=\Omega(n), \ \ n>\nu-2.
\end{equation}
According to (34), (35), (36), and (37), $\mu=\nu-1$ is the first value of $n$ for which (33) is satisfied. Then

\begin{equation*}
 {\mid}\rho_{ij}\mid \leq \mathcal{P}h^{\nu},
\end{equation*}
for any mesh sequence, where $\mathcal{P}$  is a generic constant depending on the $\{\beta_{j}\}$.\\
\newtheorem{theorem}{Theorem}[section]
\begin{theorem}
\mostafa Under the same conditions stated in Lemma 3.3, we have
\begin{equation}
{\parallel}e(x_{ij})\parallel =O(N^{-\nu}),
\end{equation}
for any mesh sequence, where
${\parallel}e(x_{ij})\parallel:={\max}\{{\mid}e(x_{ij})\mid: \ 0{\leq}i{\leq}N-1, \ 0{\leq}j{\leq}{\nu}\}.$
\end{theorem}
Proof:\\
The following might be obtained from (25) and Lemma 3.3
\begin{equation}
{\parallel}e(x_{ij}){\parallel}\leq
k_{0}\sum_{l=0}^{i-1}h_{l}{\parallel}e_{l}(x_{l}+vh_{l}){\parallel}+\mathcal{P}k_{0}h^{\nu},
\end{equation}
where
\begin{equation*}
  {\parallel}e_{l}(x_{l}+vh_{l})\parallel:={\max}\{{\mid}e_{l}(x_{l}+vh_{l})\mid:0{\leq}v{\leq}1, 0{\leq}l{\leq}i-1\}.
\end{equation*}
\begin{equation*}
  {\parallel}e(x_{ij})\parallel:={\max}\{{\mid}e(x_{ij})\mid:0{\leq}i{\leq}N-1, 0{\leq}j{\leq}\nu\}.
\end{equation*}
With the aid of Lemma 3.1, we have
\[{\parallel}e(x_{ij})\parallel=O(h^{\nu}),\]
since $h\leq N^{-1}$, then
\[{\parallel}e(x_{ij})\parallel=O(N^{-\nu}).\]
\section{Numerical Examples}
Here, we present some examples to illustrate the effectiveness and accuracy of our suggested approach.\\
{\bf Example 4.1}
\begin{equation}
y(x)=0.75x^{2}-0.5x^{2}{\ln}x+\sqrt{x}+\int_{0}^{x}\ln(x-s)y^{2}(s)ds,
\end{equation}
It is clear that this equation's exact solution is $y(x)=\sqrt{x}$.

Table 1 reports MAERR for $\nu$ = 3 based upon the present method and the technique in [23] at $x=1.0$. It is obvious that the suggested method performs similarly to the numerical solutions produced using the methodology in [23].\\
\begin{center}
\textbf{Table (1):} \ramadan The MAERR of example 4.1 for $\nu$ = 3 of the proposed technique and the approach in [23] for $x=1.0$
\end{center}
\begin{center}
\begin{tabular}{|c|c|c|}
\hline $h_{i}$&The proposed method&The method in [23] \\ \hline
0.1&1.37D-04&1.26D-03 \\ \hline
0.05&1.30D-04&3.65D-04 \\ \hline
0.025&8.59D-05&1.03D-04 \\ \hline
\end{tabular}
\end{center}
\bigskip
{\bf Example 4.2}
\begin{equation}
y(x)=0.25x^{2}-0.5x^{2}{\ln}x+x+\int_{0}^{x}\ln(x-s)(x-y(s))ds,
\end{equation}
We know the exact solution of (41) is $y(x)=x$.

Table 2, demonstrates the MAERR for $\nu$ = 3, 4, and 5 of our scheme, it is clear that as $\nu$ raises and $h$ reduces, the numerical errors of our method decay, which agrees with the theoretical results. The estimated solution and $E/h^{4}$ at $x=1.0$ using the suggested approach with $\nu = 4$, where $E$ signifies the absolute error, are numerically represented in Table (3). Table (3) demonstrates that the numerical outcomes support Theorem 3.1, which asserts that the suggested approach with $\nu = 4$ is of order four.
\section{Conclusion}
  \ \ \ \ In this research paper, an approach using a finite expansion of third-kind Chebyshev polynomials of the unknown function for second-kind VIEs with logarithmic kernels is developed. An error estimate was derived for the developed numerical procedure. Numerical experiments were given to verify the theoretical results and the efficiency of the presented scheme. This approach may be extended to solve various kinds of differential equations as well as integral and integro-differential equations.\\
\begin{center}
\textbf{Table (2):} \ramadan The proposed method's MAERR of example 4.2 for $\nu$ = 3, 4, and 5.
\end{center}
\begin{center}
\begin{tabular}{|c|c|c|c|}
\hline $h_{i}{\backslash}\nu$&3&4&5 \\ \hline
$2^{-1}$&0.77256D-02&0.41300D-03&0.10113D-03 \\ \hline
$2^{-2}$&0.39664D-02&0.17646D-03&0.39130D-04 \\ \hline
$2^{-3}$&0.15226D-02&0.56061D-04&0.11605D-04 \\ \hline
$2^{-4}$&0.51632D-03&0.15866D-04&0.31373D-05 \\ \hline
$2^{-5}$&0.16333D-03&0.42430D-05&0.81477D-06 \\ \hline
$2^{-6}$&0.49433D-04&0.11010D-05&0.20761D-06 \\ \hline
$2^{-7}$&0.14512D-04&0.28101D-06&0.52407D-07 \\ \hline
$2^{-8}$&0.41669D-05&0.71062D-07&0.13167D-07 \\ \hline
$2^{-9}$&0.11765D-05&0.17878D-07&0.33001D-08 \\ \hline
$2^{-10}$&0.32780D-06&0.44849D-08&0.82610D-09 \\ \hline
\end{tabular}
\end{center}
\bigskip
\begin{center}
\textbf{Table (3):} \ramadan The numerical values of the approximate solution and $E/h^{4}$ of the proposed method at $x=1.0$ for $\nu = 4$ of example 4.2.
\end{center}
\begin{center}
\begin{tabular}{|c|lc|}
\hline
$h_{i}$& \ \ \ \ \ \  \ \ \ \ \ \  \ \ \ \ $y(1)=1$& \\
\cline{2-3}
  &$u_{i+1}(1)=1$&$E/h^{4}$ \\ \hline
$2^{-1}$&1.0004129971&1.6520D-03 \\ \hline
$2^{-2}$&1.0001764609&2.8234D-03 \\ \hline
$2^{-3}$&1.0000560612&3.5879D-03 \\ \hline
$2^{-4}$&1.0000158658&4.0617D-03 \\ \hline
$2^{-5}$&1.0000042430&4.3449D-03 \\ \hline
$2^{-6}$&1.0000011010&4.5099D-03 \\ \hline
$2^{-7}$&1.0000002810&4.6041D-03 \\ \hline
$2^{-8}$&1.0000000711&4.6571D-03 \\ \hline
$2^{-9}$&1.0000000179&4.6866D-03 \\ \hline
$2^{-10}$&1.0000000045&4.7028D-03 \\ \hline
\end{tabular}
\end{center}


\end{document}